\documentclass[12pt]{article}

\usepackage[T1]{fontenc}
\usepackage{inputenc}
\usepackage{diagbox}
\usepackage{mathtools}
\usepackage{bbm}
\usepackage{latexsym}
\usepackage{epsfig}
\usepackage{amsmath,amsthm,amssymb,enumerate}

\usepackage[a-1b]{pdfx}
\usepackage{hyperref}

\def\wH{\tilde{H}^*_{W^*}}

\usepackage{color}

\def\cH{{\mathcal H}}

   \def\D{\Delta}
\def\e{\varepsilon}    
  \def\k{\kappa}
     
   \def\p{\pi}

\newcommand{\upp}[1]{\langle#1\rangle}

\newtheorem{theorem}{Theorem}
\newtheorem{lemma}[theorem]{Lemma}
\newtheorem{corollary}[theorem]{Corollary}
\newtheorem{Remark}{Remark}

\newcommand{\brac}[1]{\left(#1\right)}

\newcommand{\bfrac}[2]{\left(\frac{#1}{#2}\right)}

\def\cE{{\cal E}}

\newcommand{\set}[1]{\left\{#1\right\}}

\def\E{\mathbb{E}}

\def\Pr{\mathbb{P}}

\newcommand{\ignore}[1]{}

\newcommand{\THM}[2]{\
\begin{theorem}\label{#1}#2
\end{theorem}
}

\def\cA{{\mathcal A}}
\def\cB{{\mathcal B}}
\def\cC{{\mathcal C}}

\def\cE{{\mathcal E}}

\def\cH{{\mathcal H}}

\def\cS{{\mathcal S}}

\newcommand{\beq}[2]{\begin{equation}\label{#1}#2\end{equation}}

\usepackage{verbatim}
\usepackage{tikz}
\usetikzlibrary{decorations.pathmorphing}
\usetikzlibrary{positioning}
\usetikzlibrary{arrows,automata}
\usetikzlibrary{shapes.misc}
\usetikzlibrary{backgrounds}
\usetikzlibrary{arrows,shapes}

\begin{document}
\author{Tolson Bell\thanks{Research supported by NSF Graduate Research Fellowship grant DGE1745016 and DGE2140739}\\Department of Mathematical Sciences\\Carnegie Mellon University\\Pittsburgh PA 15213\and Alan Frieze\thanks{Research supported in part by NSF grant DMS1952285 }\\Department of Mathematical Sciences\\Carnegie Mellon University\\Pittsburgh PA 15213\and Trent G. Marbach\\
Department of Mathematics\\Ryerson University\\Toronto ON M5B 2K3}

\title{Rainbow Thresholds}
\maketitle

\begin{abstract}
We extend a recent breakthrough result relating expectation thresholds and actual thresholds to include some rainbow versions.
\end{abstract}
\section{Introduction}
It has long been observed that the threshold for the existence of various combinatorial objects in random graphs and hypergraphs occurs close to where the expected number of such objects tends to infinity. This informal observation has been given rigorous validation in recent breakthrough papers. First of all, Frankston, Kahn, Narayanan and Park \cite{FKNP2020} showed that under fairly general circumstances, the threshold for the existence of combinatorial objects is within a factor of $O(\log n)$ of the point where the expected number of objects begins to take off. In a follow up paper, Kahn, Narayanan and Park \cite{KNP2020} tightened their analysis for the case of the square of a Hamilton cycle and solved the existence problem up to a constant factor: a remarkable achievement, given the complexity of proofs of earlier weaker results. A key notion in this analysis is that of {\em spread}, see \eqref{spread}, first used in the paper of Alweiss, Lovett, Wu and Zhang \cite{ALWZ} that made significant progress in the resolution of the {\em Sunflower Conjecture} of Erd\H{o}s and Rado.

Spiro \cite{Spiro} describes a refinement of the notion of spread. Espuny D\'{i}az and Person \cite{Esp} generalised the approach of \cite{KNP2020} to handle some questions on spanning structures from Frieze \cite{F}.

There has been considerable research on random graphs where the edges have been randomly colored. Most notably several authors have considered the existence of rainbow colored combinatorial objects. A set of colored edges will be called {\em rainbow} if each edge has a different color. Improving on earlier results of Cooper and Frieze \cite{CF} and Frieze and Loh \cite{FL}, Ferber and Krivelevich \cite{FeKr} showed that w.h.p.~at the threshold for Hamiltonicity, randomly coloring the edges of $G_{n,p}$ with $n+o(n)$ colors yields a rainbow Hamilton cycle. Our aim in this short paper is to show that the proof in \cite{FKNP2020} can be modified to incorporate rainbow questions. We begin by summarising the results of the papers \cite{KNP2020} and \cite{FKNP2020}.

A hypergraph $\cH$ (thought of as a set of edges) is $r$-bounded if $e\in \cH$ implies that $|e|\leq r$.  The most important notion comes next. For a set $S\subseteq X=V(\cH)$ we let $\upp{S}=\set{T:\;S\subseteq T\subseteq X}$ denote the subsets of $X$ that contain $S$. Let $\upp{\cH}=\bigcup_{H\in \cH}\upp{H}$ be the collection of subsets of $X$ that contain an edge of $\cH$. We say that $\cH$ is $\k$-spread if we have the following bound on the number of edges of $\cH$ that contain a particular set $S$: 
\beq{spread}{
|\cH\cap \upp{S}|\leq \frac{|\cH|}{\k^{|S|}},\quad\forall S\subseteq X.
}
In our proof, we also use the slightly weaker definition that $\cH$ is $(K,\k)$-spread if
\beq{spread1}{
|\cH\cap \upp{S}|\leq \frac{K|\cH|}{\k^{|S|}},\quad\forall S\subseteq X.
}

Let $X_m$ denote a random $m$-subset of $X$ and $X_p$ denote a subset of $X$ where each $x\in X$ is included independently in $X_p$ with probability $p$. The following theorem is from \cite{FKNP2020}:
\THM{T2020}{
Let $\cH$ be an $r$-bounded, $\k$-spread hypergraph and let $X=V(\cH)$. There is an absolute constant $C>0$ such that if
\beq{mbound}{
m\geq\frac{(C\log r)|X|}{\k}
}
then w.h.p.~$X_m$ contains an edge of $\cH$. Here w.h.p.~assumes that $r\to\infty$.
}
\begin{Remark}\label{REM0}
Let $p=1/{\k}$ and $Z$ denote the number of edges of $\cH$ that are contained in $X_p$. Then, assuming that $\cH$ is $r$-uniform, we have from \eqref{spread} with $H\in \cH$ that $|\cH|\geq \k^r$ and then $\E(Z)=|\cH|p^r\geq 1$. This gives the connection between spread and the expected value of $Z$. Theorem \ref{T2020} inflates $p$ by a factor of order $\log r$.
\end{Remark}
Suppose now that the elements of $X$ are uniformly and independently colored from a set $Q=[q]$. A set $S\subseteq X$ is \emph{rainbow colored} if no two elements of $S$ have the same color. We modify the proof of Theorem \ref{T2020} to prove
\THM{thrainbow}{
Let $\cH$ be an $r$-bounded, $\k$-spread hypergraph and let $X=V(\cH)$ be randomly colored from $Q=[q]$ where $q\ge r$. Suppose also that (i) $\k=\Omega(r)$, that is, there exists a constant $L>0$ such that $\kappa\ge Lr$ for all valid $r$, and that (ii) $N\leq \k^2r/\log^5r$. Then there is a constant $C>0$ such that if 
\beq{mbounds}{
m\geq\frac{(C\log_2r)|X|}{\k}
}
then $X_m$ contains a rainbow colored edge of $\cH$ w.h.p. 
} 
\paragraph{Applications:}
We see that our additional restriction (i) is satisfied by many important applications and that in these cases, it implies the other extra condition, (ii). 
\begin{enumerate}[{\bf Ex. 1}]
\item Taking $X=\binom{[n]}{2}$ (the edges of $K_n$), this shows for example that if $q=n$ and $m=Kn\log n$ then w.h.p.~a randomly edge colored copy of $G_{n,m}$ contains a rainbow Hamilton cycle, see Bal and Frieze \cite{BaFr}. Here $\cH$ is the $n$-uniform hypergraph with $(n-1)!/2$ edges, one for each Hamilton cycle of $K_n$, $r=n$ and $\k\geq n/e$. 

\item Dudek, English and Frieze \cite{DEF} studied rainbow Hamilton cycles in random hypergraphs. Theorem \ref{thrainbow} improves Theorem 6 of that paper. Here $X=\binom{[n]}{k}$ and $\cH$ is the $n/(k-1)$-uniform hypergraph with $\frac{(k-1)n!}{2n(k-2)!^{n/(k-1)}}$ edges, one for each loose Hamilton cycle. We can take $\k=\Omega(n^{k-1})$. 

\item Similarly, if $q=n$ and $m=C_1n\log n$ and $T$ is an $n$-vertex tree with bounded maximum degree $\D=O(1)$ then w.h.p.~$G_{n,m}$ contains a rainbow copy of $T$. It suffices to take $\k=n/\D$. The uncolored version is due to Montgomery \cite{M1}. One can easily extend this result to hypergraphs. Spanning trees can be replaced by spanning {\em cacti/hypertrees}. Let $K_{n,k}$ denote the complete $k$-uniform hypergraph on vertex set $[n]$. A single edge is a cactus and a cactus $C_1$ with $m+1$ edges is obtained from a cactus $C_0$ with $m$ edges by selecting a vertex $v\in V(C_0)$ and adding an edge $\set{v=v_1,v_2,\ldots,v_k}$ where $v_2,v_2,\ldots,v_k\notin V(C_0)$. A cactus is spanning if its vertex set $V(C_n)=[n]$. A cactus with $m$ edges contains $m(k-1)+1$ vertices and so we need $(k-1)\mid (n-1)$. Let $C_n$ be a sequence of spanning cacti of $K_{n,k}$ all of maximum degree $\D=O(1)$. We take $X=\binom{[n]}{k}$ and the edges of $\cH$ correspond to the copies of $C_n$ in $K_{n,k}$. We prove that \eqref{spread} holds with $\k=\binom{n-1}{k-1}/\D$. If $S\subseteq X$ is not isomorphic to a subset of $E(C_n)$ then $E(\cH)\cap \upp{S}=\emptyset$ and \eqref{spread} holds. Suppose then that $S\subseteq X$ is isomorphic to a subset of $E(C_n)$. Then, where $\p$ is a random permutation of $[n]$ and 
\[
\p(C_n)=([n],\set{\set{\p(v_1),\p(v_2),\ldots,\p(v_k)}:\set{v_1,v_2,\ldots,v_k}\in E(C_n)},
\]
we have
\beq{bddtree}{
\frac{|E(\cH)\cap \upp{S}|}{|E(\cH)|}=\Pr(S\subseteq \p(C_n))\leq \k^{-|S|}.
}
To see \eqref{bddtree}, suppose that $S=\set{e_1,\ldots,e_s}$. Fix a vertex $v$ that occurs in one of the edges of $S$ and suppose that $v\in e_1,e_2,\ldots,e_\ell$. Then we have 
\[
\Pr(e_1,e_2,\ldots,e_\ell\in \p(C))\leq \frac{(\D)_{\ell}}{\brac{\binom{n-1}{k-1}}_\ell}\leq \k^{-\ell}.
\]
Furthermore, after removing $e_1,e_2,\ldots,e_\ell$ and $v$ from the cactus $C$ we can apply induction. Applying \eqref{mbounds} we see that $O(n\log n)$ randomly colored random edges are needed for a rainbow coloring. 

\item We also obtain a rainbow version of Shamir's problem, Corollary 1.2 of Bal and Frieze \cite{BaFr}. Here $X=\binom{[n]}{k}$ and $\cH$ is the $n/k$-uniform hypergraph with $\frac{n!}{(n/k)!k!^{n/k}}$ edges, one for each perfect matching. It suffices to take $\k=n^{k-1}/k!$. 
\end{enumerate}
The paper \cite{KNP2020} focusses exclusively on the square of Hamilton cycles. It removed a $\log n$ factor. So, given this we know that there are constants $0<c_1<c_2$ such that w.h.p.~$G_{n,c_1n^{3/2}}$ does not contain the square of a Hamilton cycle and that  w.h.p.~$G_{n,c_2n^{3/2}}$ does contain the square of a Hamilton cycle. The proof is similar to that of Theorem \ref{T2020}. Here $\k=O(n^{1/2})$ and the constraint that $\k=\Omega(r)$ prevents us from generalising this result in exactly the same way. On the other hand, Bell and Frieze \cite{BellF} prove a rainbow version using a different but related argument.

\begin{Remark}\label{rem00}
An obvious line of attack here is to (i) replace each edge $e=\set{v_1,v_2,\ldots,v_k}$ from the original hypergraph $\cH$ by a set of edges $\{\set{(v_1,c_1),(v_2,c_2),\ldots,(v_k,c_k)}\subseteq X^*=X\times [q]$ with all $c_i$ distinct$\}$, (ii) verify that the new hypergraph is $\k'$-spread for some value $\k'$ and then apply  \cite{KNP2020} or \cite{FKNP2020}. The appropriate value of $\k'$ seems to be $q\k/e$, see \eqref{starspread}. The only problem with this approach is that a random $m$-subset of $X^*$ may not correspond to a randomly colored subset of $X$. Let $\cB$ denote the subsets of $X_m^*$ that contain a pair $(x,c_1),(x,c_2)$ i.e. where $x\in X$ has been given two colors. If $N=|X|$ then the expected number of such pairs in a random subset $X_m^*$ is $\approx \frac{Nq^2}2\cdot\bfrac{m}{qN}^2=O\bfrac{N\log^2r}{\k^2}$. So we see immediately that Theorem \ref{T2020} holds if $\frac{\k}{N^{1/2}\log r}\to\infty$. Otherwise, we have to work around the problem. 
\end{Remark} 
\section{Proof of Theorem \ref{thrainbow}}\label{SML}
We use a superscript $*$ to indicate colored objects. Let $X^*=X\times [q]$. For $x^*=(x,c)\in X^*$, define $c(x^*)=c$ (called the {\em color} of $x^*$) and $\xi(x^*)=x$. For $S^*\subseteq X^*$, define $\xi(S^*)=\set{\xi(x):x\in S^*}$. We say that a set $S^*\subseteq X^*$ is {\em rainbow} if $c(x^*)\neq c(y^*)$ for all distinct $x^*,y^*\in S^*$. Let 
\[
\cE^*=\set{W^*\subseteq X^*:\xi(x^*)\ne\xi(y^*)\text{ for any distinct }x^*,y^*\in W^*}.
\]
We say that $A^*,B^*\subseteq X^*$ are {\em compatible}, and write $A^*\sim B^*$, if $A^*\cup B^*\in\cE^*$, that is, for all $a^*,b^*\in A^*\cup B^*$, $\xi(a^*)= \xi(b^*)$ implies $c(a^*)=c(b^*)$.

\subsection{Single-stage Refinement}\label{stage}
In this section, following the general framework of \cite{ALWZ,FKNP2020,PP2022}, our main result is Lemma \ref{ML}, which will later be used iteratively to bound the number of sets of $\cH$ in Theorem \ref{thrainbow} that are ``bad''.

Suppose now that $\cH$ is an $r$-uniform $(K,\k)$-spread hypergraph and that each $H\in\cH$ can be colored from one of a set of colors $Q_H\subseteq Q$ where $|Q_H|=q$. Let 
\beq{Hstar}{
\cH^*=\set{(H,c_H):H\in\cH,c_H:H\to Q_H,c_H\text{ is 1-1}}.
}
Thus $\cH^*$ consists of rainbow colored edges.
\begin{lemma}\label{qk}
With $\cH^*$ as defined in \eqref{Hstar}, $\cH^*$ is $(K,q\k/e)$-spread.
\end{lemma}
\begin{proof}
Let $X^*=X\times Q$. If $S^*\subseteq X^*$ is rainbow and $|S^*|=s$, then as $\cH$ is $\kappa$-spread,
\beq{starspread}{
|\cH^*\cap\upp{S^*}|=\sum_{\substack{H\in \cH\\H\supseteq \xi(S^*)}}(q-s)_{r-s}\leq  (q-s)_{r-s}\frac{|\cH|}{\k^s}= \frac{(q-s)_{r-s}}{(q)_r}\frac{|\cH^*|}{\k^s}\leq \frac{e^s|\cH^*|}{(q\k)^s}.
}
If $S^*$ is not rainbow or $|\xi(S^*)|\neq |S^*|$ then $\cH^*\cap\upp{S^*}=\emptyset$. Thus, $\cH^*$ is $q\k/e$ spread. 
\end{proof}
For the remainder of this section, let $\cH^*$ be any multi-hypergraph (that is, allowing repeated edges) on $X^*=X\times[q]$ that is $r$-bounded and $(K,q\kappa/e)$-spread for some $r\le q,\k>1$ and $K\geq 1$. Assume also that every $H^*\in\cH^*$ is in $\cE^*$ and rainbow, that is, $\forall~(x_1,c_1),(x_2,c_2)\in H^*$, $x_1=x_2\iff c_1=c_2$. 

Set $p=\frac{C_0}{\kappa}$ for $C_0$ large, $N=|X|,$ and let $m=Np$. Let $W^*_m$ be chosen randomly from $\cE^*_m=\set{E^*\in \cE^*:|E^*|=m}$. In other words, $W^*_m$ is the same as choosing $m$ random vertices and randomly coloring them. Let
\beq{JBx}{
\cH^*_{W^*}=\set{H^*\in \cH^*:H^*\sim W^*},
}
that is, the set of rainbow hyperedges that are compatible with $W^*$.

For $H^*\in\cH^*$ such that $H^*$ is compatible to $W^*$, let $T^*=T^*(H^*,W^*)$ be $G^*\setminus W^*$ for a set $G^*\in\cH^*$ such that $G^*\subseteq H^*\cup W^*$ and $|G^*\setminus W^*|$ is minimized. For $1\le t\le r$, let $\cC_t(W^*)=\{H^*\in\cH^*:\mbox{$H^*\sim W^*,|T^*(H^*,W^*)|=t$}\}$.

\begin{lemma}\label{cCt}
For $\cH^*$ as above, $W^*$ chosen uniformly from $\cE^*_m$, and $1\le t\le r$,
\beq{ESWBad}{
\E\brac{|\cC_t(W^*)|}\le \frac{K2^re^t}{C_0^{t}}|\cH^*|}
\end{lemma}
\begin{proof}
Note that there are $\binom Nmq^m$ equally likely choices for $W^*$. Thus, it suffices to show that 
\[
\sum_{W^*\in\cE^*_m}|\cC_t(W^*)|\le\binom Nmq^m\frac{K2^re^t}{C_0^{t}}|\cH^*|.
\]
 We will give a procedure that uniquely specifies every possible $(W^*,H^*)$ pair such that $|T^*|=t$. Fix some function $\chi:2^{X^*}\rightarrow\cH^*$ such that $\chi(Y^*)\subseteq Y^*$ whenever possible.

\begin{enumerate}[\textit{Step} 1.]
\item Specify $\xi(W^*\cup T^*(H^*,W^*))$. There are $\binom{N}{m+t}\le\binom Nm\bfrac{\k}{C_0}^{t}$ choices here.

\item Specify $Z^*=W^*\cup T^*$. There are $q^{m+t}$ choices here.

\item Note that we must have $T^*\subseteq\chi(Z^*)$ by the minimality of $T^*$, as $\chi(Z^*)$ is a valid choice for $G^*$, $\chi(Z^*)\backslash W^*\subseteq T^*=T^*\setminus W^*$, and $|T^*\setminus W^*|\le|\chi(Z^*)\backslash W^*|$. Thus, we can specify $T^*(H^*,W^*)\subseteq\chi(Z^*)$ with at most $\binom{|\chi(Z^*)|}{t}\le 2^{|\chi(Z^*)|}\le 2^r$ choices.

\item Specify $H^*\supseteq T^*$. As $\cH^*$ is $(K,q\k/e)$-spread and $|T^*|=t$, there are at most $(q\k/e)^{-t}K|\cH^*|$ choices
\end{enumerate}
Thus, $(W^*,T^*)$ has been specified with at most 
\[
\binom Nm\bfrac{\k}{C_0}^{t}q^{m+t}2^r \bfrac{e}{q\k}^{t}K|\cH^*|=\binom Nmq^m\frac{K2^re^t}{C_0^{t}}|\cH^*|
\]
 choices.
\end{proof}

\begin{lemma}\label{ML}
Let $\cH^*$, $W^*$ be as above. Now, let $H^*\in\cH^*$ be good with respect to $W^*$ if $H^*\sim W^*$ and $|T^*(W^*,H^*)|<r/2$. Let 
\[
\text{$\cS$ be the event $|\set{H^*\in \cH^*:H^*\text{ is good}}|<(1-\e)|\cH^*|(1-p)^r$.}
\]
Then
\[
\Pr(\cS)\leq \exp\set{-\frac{\e^2\k^2q(1-p)}{16e^3N}} +\frac{2K}{C_0^{r/3}\e\brac{1-p}^r}
\]
\end{lemma}
\begin{proof}
Let $\wH=\set{H^*\in \cH^*:\xi(H^*)\cap \xi(W^*)=\emptyset}$. Let $\cS_1$ be the event that $|\wH|\leq (1-\e/2)|\cH^*|(1-p)^r$ and let $\cS_2$ be the event that $\sum_{t>r/2}^r|\cC_t(W^*)|\geq \e|\cH^*|(1-p)^r/2$. Then $\cS\subseteq \cS_1\cup \cS_2$. This is because every $H^*\in \wH$ is compatible with $W^*$ and if $\cS_1$ fails then there must be at least $\e|\cH^*|(1-p)^r/2$ sets $H^*$ with $|T^*(W^*,H^*)|\geq r/2$

We will use the Janson inequality. If $Z=|\wH|$ then we have 
\[
\E(Z)=|\cH^*|(1-p)^r.
\]
If $H_1^*,H_2^*\in \cH^*$ then we write $H_1^*\sim H_2^*$ if $\xi(H_1^*)\cap \xi(H_2^*)\neq\emptyset$. Then,
\begin{align*}
\D&=\sum_{\substack{H_1^*,H_2^*\in\cH^*\\H_1^*\sim H_2^*}}\Pr(H_1^*,H_2^*\in\wH)\\
&\leq \sum_{s=1}^r\sum_{\substack{S^*\subseteq X^*\\|S^*|=s}}\sum_{H_1^*,H_2^*\supseteq S^*}\Pr(H_1^*,H_2^*\in\wH)\\
&\leq\sum_{s=1}^r\binom{qN}{s}\bfrac{|\cH^*|}{(q\k/e)^s }^{2}(1-p)^{2r-s}\\
&\leq \E(Z)^2\sum_{s=1}^r\bfrac{qNe}{s}^s\bfrac{e^2}{q^2\k^2 (1-p)}^{s}\\
&= \E(Z)^2\sum_{s=1}^r\bfrac{Ne^3}{q\k^2s(1-p)}^s\\
&\leq \E(Z)^2\bfrac{2Ne^3}{q\k^2(1-p)}.
\end{align*}
And by Janson's inequality,
\[
\Pr\brac{Z\leq \brac{1-\tfrac{\e}{2}}\E(Z)}\leq \exp\set{-\frac{\e^2\E(Z)^2}{8\D}} \leq \exp\set{-\frac{\e^2\k^2q(1-p)}{16e^3N}}.
\]
As for $\cS_2$, by Lemma \ref{cCt}, 
\[
\E\brac{\sum_{t>r/2}^r|\cC_t|}\le K|\cH^*|\sum_{t>r/2}^r\frac{2^r}{C_0^{t}}\leq \frac{K}{C_0^{r/3}}|\cH^*|.
\]
And thus by Markov's inequality,
\[
\Pr(\cS_2)\le \frac{2K}{C_0^{r/3}\e(1-p)^r}.
\]
\end{proof}
We next let $\cH^*_{W^*}=\set{T(W^*,H^*):\text{$H^*$ is good with respect to $W^*$}}$. Then,
\begin{lemma}\label{spreadS}
If $\cH^*$ is $(K,q\k/e)$ spread and $\cS$ does not occur then $\cH_{W^*}^*$ is $\brac{K/((1-\e)(1-p)^r),q\k/e}$ spread.
\end{lemma}
\begin{proof}
\[
|\cH^*_{W^*}\cap\upp{S^*}|\leq |\cH^*\cap\upp{S^*}|\leq \frac{K|\cH^*|}{(q\k/e)^s}\leq \frac{K|\cH^*_{W^*}|}{(1-\e)(1-p)^r(q\k/e)^s}.
\]
\end{proof}
\subsection{Completing the proof}\label{SMT}
Let $\cH$ be as in Theorem \ref{thrainbow}, an $r$-bounded, $\kappa$-spread hypergraph on a set $X$ of size $N$ such that $r\le C_1\k$ for some constant $C_1$. Then let 
\beq{hcstar}{
\cH^*=\set{S^*\subseteq X^*:\xi(S^*)\in\cH,\text{ and }c(x^*)\ne c(y^*)\text{ and }\xi(x^*)\ne\xi(y^*)\text{ for any distinct }x^*,y^*\in S^*},
}
be the set of possible rainbow edges of $\cH$.

We follow the basic idea of \cite{FKNP2020} and \cite{PP2022} and account for the coloring. We choose a small randomly colored random set $W^*$, throw away the $H^*$ that are not compatible with $W^*$, and replace the other $H^*$ with $T^*(W^*,H^*)$. We argue that the way we do this causes most $H^*$ to become relatively smaller. We then repeat the argument with respect to the hypergraph consisting of these $T^*$. In this way, we build up $W^*$ piece by piece and find members of $\cH^*$ which are mostly contained in $W^*$. After $O(\log r)$ iterations, one of those pieces is likely to be fully contained in a larger $W^*$. It is important to realise that the edge sets of the hypergraphs encountered in these iterations are multi-sets, i.e., the same edge can be repeated many times. Recall that asymptotics refer to $r$.

Let $\ell_0$ be the least $i$ such that $r/2^i<\log r$ and let $\ell$ be the least $i>\ell_0$ such that $r/2^i<1$. Note that $r_{\ell}\geq1/2$. Let $p_i=C_i/\kappa$ where 
\begin{align*}
C_i&=\begin{cases}C_0&i\leq \ell_0\\\frac{\log r}{\log\log r}&\ell_0<i\leq \ell\end{cases};\qquad \e_i=\begin{cases}\frac{1}{(i+1)^2}&1\leq i\leq \ell_0.\\\frac{1}{(\log\log r)^2}&i>\ell_0.\end{cases};\qquad r_i=\frac{r}{2^i} \text{ for }i\geq 0.\\
K_i&=K\prod_{j=1}^i\frac{1}{(1-\e_j)\brac{1-p_j}^{r_j}}.
\end{align*}
Here $C_0$ is some sufficiently large constant. Note that 
\[
\sum_{j=1}^{\ell} r_jp_j\leq \frac{2C_0r}{\k}+\sum_{j=\ell_0+1}^{\ell}\frac{\log r}{\k\log\log r}=\frac{2C_0r}{\k}+o(1)\text{ and }\prod_{i=\ell_0+1}^{\ell}\frac{1}{1-\frac{1}{(\log\log r)^2}}=e^{o(1)}.
\]
So, using $1/(1-x)\leq e^{2x}$ for small $x$, we obtain
\beq{KK}{
1\leq K_i\leq \exp\set{2\sum_{t=1}^\infty\frac{1}{(t+1)^2}+o(1)+2\sum_{j=1}^{\ell} r_jp_j}\leq A,
}
where $A=B+e^{3C_0/L}$ for some absolute constant $B$. (We have used the fact that $\k$-spread can be interpreted as $(1,\k)$-spread.)

Let $X_0=X$ and $X_0^*=X\times [q]$. For $i=1,\dots ,\ell$, let $X_i=X_{i-1}\setminus \xi(W_i^*)$, where $W^*_i$ is chosen uniformly at random from $\cE^*_{Np_i}$ on $X_{i-1}^*$ (that is, $\set{W^*_i\in \cE^*_{Np_i}:W_i^*\subseteq X_{i-1}^*}$). Let
\[
\cH_i^*=\{T^*(H^*,W^*_i):\mbox{$H^*\in \cH_{i-1}^*$, $H^*\sim W^*_i$, $|T^*|\le r_i$}\}.
\]
Thus $\cH_i^*$ is an $r_i$-bounded collection of subsets of $X_i$. 

Note that $\exists H^*\in \cH^:H^*\subseteq W_i^*$ whenever some $T^*(H^*,W_i^*)=\emptyset$ for $H^*\in\cH_{i-1}^*$. As $r_{\ell}<1$, eventually we will either have some $T^*(H^*,W_i^*)=\emptyset$ or $|\cH_i^*|=0$. So to prove the claim, we just need to show
\beq{toshow'}{
\Pr(\cH^*_{\ell}=\emptyset)=o(1),
}
(where the $\Pr$ refers to the entire sequence $W^*_1\dots W^*_{\ell}$).

Let $\cA_i$ be the event that $\cS_i$ does not occur, where $\cS_i$ is the $\cH^*_i$ equivalent of the event $\cS$ in Lemma \ref{ML} and let $\cA_{\leq i}=\bigcap_{j=0}^i\cA_ j$. 
\beq{fA}{
\text{If $\cA_{\leq \ell}$ occurs then }|\cH^*_{\ell}|\geq |\cH^*|\prod_{i=1}^{\ell}(1-\e_i)(1-p_i)^{r_i}>0.
}
Lemma \ref{ML} implies that
\begin{align*}
\Pr(\cA_{\leq i}\mid \cA_{\leq i-1})&\geq 1-\exp\set{-\frac{\e_i^2\k^2q(1-p)}{16e^3N}}-\frac{2K_i}{C_i^{r_i/3}\e_i\brac{1-\frac{C_i}{\k}}^{r_i}}\\
&\geq 1-\exp\set{-\frac{\k^2q}{400(i+1)^4N}}-\frac{2A}{C_i^{r_i/3}\e_i\brac{1-\frac{C_0}{\k}}^{r}}>0\\
&\geq 1-r^{-1/400}-\frac{2A}{C_i^{r_i/3}\e_i\brac{1-\frac{C_0}{\k}}^{r}}>0.
\end{align*}
Note that $C_i^{r_i/3}\e_i\geq C_0^{\frac13\log r}/\log^2r=\Omega(r)$ for $i\leq \ell_0$ and $C_i^{r_i/3}\e_i\geq (\log r/\log\log r)^{1/3}/(\log\log r)^2=\Omega(\log^{1/4}r)$ for $\ell_0<i\leq \ell$.

So,
\begin{align*}
\Pr(\cA_{\leq\ell})&\geq \prod_{i=1}^{\ell_0}\brac{1-r^{-1/400}-O\bfrac{1}{r}} \prod_{i=\ell_0+1}^\ell\brac{1-\frac{e^{-C_0/L}}{(i+1)^2(1-e^{-C_0/L})}-O\bfrac{1}{\log^{1/4}r}}\\
&= (1-o(1)).
\end{align*}
Applying \eqref{toshow'} and \eqref{fA} completes the proof of Theorem \ref{thrainbow}.
%%%%%%%%%%%%%%%%%%%%%%%%%%%%%%%%%%%%%%%%%%%%%%%%%%%%%%%%%%%%%%%%%%%%%%%%%%%%%%%
\section{Final Remarks}
An earlier version claimed to have full generalisations of the threshold results of \cite{FKNP2020}, \cite{KNP2020}. Excellent reviewing pointed out significant errors and we can only claim a partial generalisation of \cite{FKNP2020} due to the restriction $\k=\Omega(r)$ and $N\leq \k^2r/\log^5r$. We nevertheless are confident that the threshold results of these papers can be generalised to rainbow versions.

{\bf Acknowledgement} We thank Sam Spiro and Lutz Warnke and Erlang Surya for pointing out issues in previous versions. We also thank the referees for their help and dedication in correcting the paper.

\section{Subsequent developments}
Subsequent to the initial production of this paper, there have been two significant developments.

First of all Han and Yuan \cite{HY} have proved a version of Theorem \ref{thrainbow} without the extra conditions. This does not lead to any significant new applications.

Bell and Frieze \cite{BellF} amend the proof in \cite{KNP2020} to prove
\begin{theorem}
Let $\e,\e_1>0$ be arbitrary positive constants. Suppose that $\cH$ is a $\k$-spread, $r$-uniform and edge transitive hypergraph, subject to a technical condition $\cC$. Let $X=V(\cH)$ be randomly colored from $Q=[q]$ where $q\geq (1+\e_1)r$. Then there exists $C=C(\e,\e_1)$ {such that for sufficiently large $r,\k$,}
\beq{strictly}{
m\geq \frac{C|X|}{\k} \text{ implies that }{\Pr(X_{m}^*\text{ contains a rainbow colored edge of $\cH$}})\geq 1-\e.
}
\end{theorem}
Powers of Hamilton cycles give rise to hypergraphs satisfying $\cC$ and so $p=n^{-1/k}$ is the threshold for the $k$th power of a Hamilton cycle, given $q=(1+\e_1)n$ random colors.
\end{document}